\documentclass[12pt,oneside]{amsart}
\usepackage{amsthm,amsfonts,amssymb,euscript}

\renewcommand{\emptyset}{\varnothing}
\newcommand{\mt}[1]{\operatorname{#1}}
\newcommand{\red}{\operatorname{red}}
\newcommand{\Supp}{\operatorname{Supp}}
\newcommand{\Diff}{\operatorname{Diff}}
\newcommand{\down}[1]{\left\lfloor #1\right\rfloor}

\newcommand{\discr}{\operatorname{discr}}
\newcommand{\totaldiscr}{\operatorname{totaldiscr}}
\newcommand{\CC}{{\mathbb C}}
\newcommand{\ZZ}{{\mathbb Z}}
\newcommand{\QQ}{{\mathbb Q}}
\newcommand{\PP}{{\mathbb P}}
\newcommand{\NN}{{\mathbb N}}

\newcommand{\OOO}{{\mathcal O}}
\newcommand{\EEE}{{\EuScript{E}}}
\newcommand{\Msm}{{\Phi}_{\mt{\mathbf{sm}}}}
\newcommand{\qq}{\mathbin{\sim_{\scriptscriptstyle{\mathbb Q}}}}

\newtheorem{theorem}[subsection]{Theorem}
\newtheorem{conjecture}[subsection]{Conjecture}
\newtheorem{proposition}[subsection]{Proposition}
\newtheorem{lemma}[subsection]{Lemma}
\newtheorem{corollary}[subsection]{Corollary}
\theoremstyle{definition}

\newtheorem{example}[subsection]{Example}
\newtheorem{proposition-definition}[subsection]{Proposition-Definition}
\theoremstyle{remark} 
\newtheorem{remark}[subsection]{Remark}

\title{A note on log canonical thresholds}
\author{Yu.~G.~Prokhorov}
\thanks{This work was partially supported by the grant
INTAS-OPEN-97-2072}
\address{Department of Mathematics
(Algebra Section), Moscow Lomonosov University, 117234 Moscow,
Russia}
\email{prokhoro@mech.math.msu.su}

\address{Current address: Max-Plank-Institut f\"ur Mathematik
Vivatsgasse 7, 53111 Bonn, Germany}
\email{prokhoro@mpim-bonn.mpg.de}

\begin{document}
\begin{abstract}
We prove that the largest accumulation point of the set
$\mathcal{T}_3$ of all three-dimensional  
log canonical thresholds $c(X,F)$ is $5/6$.
\end{abstract}
\maketitle
\section{Introduction}
Let $(X,\Omega)$ be a log variety and let $F$ be an effective
non-zero Weil
$\QQ$-Cartier divisor on $X$. Assume that $(X,\Omega)$ has at
worst log canonical singularities. \textit{The log canonical
threshold} of $F$ with respect to $(X,\Omega)$ is defined by
\[
c(X,\Omega,F)=\sup \left\{c \mid \text{$(X,\Omega+cF)$
is log canonical}\right\}.
\]
It is known that $c(X,\Omega,F)$ is a rational number 
from the interval $[0,1]$ (see \cite{Ko}). 
We frequently write $c(X,F)$ instead of $c(X,0,F)$.

For each $d\in \NN$ define the set $\mathcal{T}_d\subset [0,1]$ by
\[
\mathcal{T}_d:= \left\{c(X,F)\ \left|\ \begin{array}{c}
\text{$\dim X=d$, $X$ has only log canonical singularities}\\ 
\text{and $F$
is an effective non-zero Weil $\QQ$-Cartier divisor}
\end{array}\right.\right\}.
\]
The structure of $\mathcal{T}_d$ is interesting for applications 
to the problem of termination some inductive procedures appearing 
in the Minimal Model Program \cite{Sh}, \cite{Ut}.   
The interest in log  canonical thresholds was also inspired
in connection with the complex singular index and Bernstein-Sato
polynomials (see \cite{Ko}).

\begin{conjecture}[{\cite{Sh}}]
\label{conj}
$\mathcal{T}_d$ satisfies the ascending chain condition, i.e. 
any increasing chain of elements terminates.
\end{conjecture}

The set $\mathcal{T}_2$ is completely described (see 
\cite{Ku1}). Concerning $\mathcal{T}_3$ it is known 
the following:
\begin{enumerate}
\item
Conjecture \ref{conj} holds true for $\mathcal{T}_3$ \cite{A0},
\cite[Ch. 18]{Ut};
\item
$\mathcal{T}_3\cap (41/42,1)=\emptyset$ \cite{Ko1};
\item
$\mathcal{T}_3\cap [6/7,1]$ is finite \cite{Pr2}.
\end{enumerate} 

Actually, the structure of $\mathcal{T}_d$ is rather complicated:
it has a lot of accumulation points  \cite[8.21]{Ko}. 
However adopting Conjecture \ref{conj} we see that  $\mathcal{T}_d$
is discrete near $1$.  

Our main result is the following theorem which generalizes 
the result of \cite{Pr2}. 

\begin{theorem}
\label{main}
The largest accumulation value of $\mathcal{T}_3$ is $5/6$.
\end{theorem}

\begin{remark} (i)
The two-dimensional analog of our theorem easily follows from 
the description of $\mathcal{T}_2$ (\cite{Ku1}):
the largest accumulation value of $\mathcal{T}_2$ is $1/2$. 

(ii) T. Kuwata 
described the set of all values  $c(\CC^3,F)$ in the interval 
$[5/6,1]$, where $F$ is a 
hypersurface in $\CC^3$. His proof is done by studying the local 
equation of $F$. Our proof uses quite different method and 
based on Alexeev's result \cite{A}.  
\end{remark}

The essential part of the proof is to show the finitedness of 
$\mathcal{T}_3\cap [5/6+\epsilon,1]$
for any $\epsilon>0$.
The  easy example below  shows that $5/6$ is an accumulation point of 
$\mathcal{T}_3$. 

\begin{example}
\label{ex-1}
Let $X=\CC^3$ and let $F_r$ be the hypersurface given by
$x^2+y^3+z^r$, $r\ge 7$. This singularity is quasihomogeneous.  
By \cite[8.14]{Ko} we have
$c(\CC^3,F_r)=5/6+1/r$. Thus $\lim_{r\to \infty} c(\CC^3,F_r)=5/6$. 
\end{example}

\par\medskip\noindent
\textsc{Acknowledgments.}
This work was completed during my stay at 
Max-Planck-Institut f\"ur Mathematik.
I would like to thank MPIM for hospitality and support.
I am grateful to Dr. O. Fujino for pointing out 
several inaccuracies in the first draft. 

\section{Preliminary results}
All varieties are assumed to be algebraic varieties defined over
the field $\CC$. A \textit{log variety} (or a \textit{log pair}) $(X,D)$
is a normal quasiprojective variety $X$ equipped with a \textit{boundary},
a $\QQ$-divisor $D=\sum d_iD_i$ such that $0\le d_i\le 1$ for all $i$.
We use terminology, definitions and abbreviations of the Minimal Model 
Program \cite{Ut}. 

\begin{proposition-definition}[{\cite[\S 3]{Sh}, \cite[Ch. 16]{Ut}}]
Let $(X,S+B)$ be a log variety, where 
$S=\down{S+B}\ne 0$ and divisors $S$, $B$ have no common components. 
Assume that
$K_X+S$ is lc in codimension two. Then there is a naturally defined 
effective $\QQ$-divisor $\Diff_S(B)$ on $S$ called the 
\emph{different} of $B$ such that 
$$
K_S+\Diff_S(B)\qq (K_X+S+B)|_S.
$$
\end{proposition-definition}

\subsection{}
Let $\Phi$ be a subset of $\QQ$. For a $\QQ$-divisor $D=\sum
d_iD_i$, we write $D\in \Phi$ if $d_i\in\Phi$ for all $i$.
Define the following sets
\[
\begin{array}{l}
\Msm:=\bigl\{1-1/m \mid m\in\NN\cup\{\infty\}\bigr\},\\
\Msm^\alpha:=\Msm\cup [\alpha,1],\quad\text{for $\alpha\in[0,1]$}.
\end{array}
\]
We distinguish them because they 
are closed under some important operations (see e.g. 
Corollary \ref{Msma}
below).
Usually the numbers from $\Msm$ are called 
\emph{standard}.

\begin{proposition}[{\cite[Prop. 3.9]{Sh}}]
\label{iso}
Let $(X,S)$ be a $d$-dimensional plt log variety, where $S$ is 
integral. Let $W\subset S$ be an irreducible subvariety of
codimension $1$. Then near the general point $P\in W$ there is an
analytic isomorphism
\begin{multline}
\label{m}\qquad 
(X,S,W)\simeq
\Bigl((\CC^d,\{x_1=0\},\{x_1=x_2=0\})/\ZZ_m(1,q,0\dots,0)\Bigr),\\
\text{where}\quad m,\ q\in\NN,\ \gcd(m,q)=1.
\end{multline}
\end{proposition}

\begin{corollary}[{\cite[3.10, 3.11]{Sh}}]
\label{coeff_Diff_1}
Let $(X,S+B)$ be a log variety, where $S:=\down{S+B}$
and divisors $S$, $B$ have no common components.
Assume that $(X,S)$ is plt. Let $W\subset S$ be an irreducible
subvariety of codimension $1$. If $B=\sum b_iB_i$, then the
coefficient of $\Diff_S(B)$ along $W$ is equal to
\begin{equation}
\label{coeff_Diff}
1-\frac1m+\sum_{B_i\supset W} \frac {n_ib_i}{m},
\end{equation}
where $m$ is such as in \eqref{m} and $n_i\in\NN$. Moreover, if
$(X,S+B)$ is plt and $B\in [1/2,1]$, then 
there is at most one component $B_i$ of $B$ containing $W$ and $n_i=1$.
\end{corollary}

\begin{corollary}[{\cite[3.11, 4.2]{Sh}}]
\label{Msma}
Let $(X,S+B)$ be a log variety, where $S:=\down{S+B}$
and divisors $S$, $B$ have no common components.
Assume that $(X,S)$ is plt and $(X,S+B)$ is plt.
Take $\alpha\in [0,1]$.
If $B\in\Msm^\alpha$, then $\Diff_S(B)\in \Msm^\alpha$.
\end{corollary}

\begin{proposition-definition}[{\cite{Pr1}}]
\label{constr-plt-blowup}
Let $(X,D)$ be a log variety such that 
$(X,D)$ is lc but not plt, $X$ is klt and
$\QQ$-factorial. Assume the log MMP
in dimension $\dim(X)$. Then there exists a blow-up $f\colon Y\to
X$ such that
\begin{enumerate}
\item
the exceptional set of $f$ contains an unique prime divisor $S$;
\item
$K_Y+D_Y=f^*(K_X+D)$ is lc, where $D_Y$ is the proper transform of
$D$;
\item
$K_Y+S+(1-\varepsilon)D_Y$ is plt and anti-ample over 
$X$ for any $\varepsilon>0$;
\item
$Y$ is $\QQ$-factorial and $\rho (Y/X)=1$.
\end{enumerate}
Such a blow-up we call an \textit{inductive blow-up} of $(X,D)$.
\end{proposition-definition}

\section{Lemmas}

\begin{lemma}
\label{P1}
Let $\Lambda$ be a boundary on $\PP^1$ such that 
$\Lambda\in\Msm^{5/6}$ and $K_{\PP^1}+\Lambda\equiv 0$. 
Then $\Lambda\in \Msm\cap [0,5/6]\cup \{1\}$.
\end{lemma}
\begin{proof}
Write $\Lambda=\sum \lambda_i\Lambda_i$. Then $\lambda_i\in \Msm^{5/6}$
and $\sum \lambda_i=2$.
If $\down{\Lambda}\neq 0$, then there are only two possibilities:
$\lambda_1=\lambda_2=1$ and $\lambda_1=2\lambda_2=2\lambda_3=1$.
Otherwise $\lambda_i<1$ and
easy computations give us 
$\lambda_i\le 5/6$, so $\lambda\in\Msm$.
\end{proof}

\begin{lemma}
\label{LCT-2}
Let $\left(S,\Delta =\sum \delta_{i}\Delta_{i}\right)$ be a lc log
surface such that $\delta_{i}\in \Msm^{5/6}$ and let $C$ be an
effective Weil divisor on S. Then either $c(S,\Delta,C)\leq 5/6$
or $c(S,\Delta,C)=1$.
\end{lemma}

\begin{proof}
Put $c:=c(S,\Delta,C)$. Assume that $5/6<c<1$. By \cite[8.5]{Ko}
there is an exceptional divisor $E$ such that $a(E, \Delta
+cC)=-1$ and $a(E, \Delta)>-1$. 
Put $P:=\mt{Center}(E)$. Regard $S$ as a germ near $P$.

Let $\varphi
\colon \widetilde S\to S$ be an inductive blowup of $(S,\Delta
+cC)$. Write
\begin{equation*}
K_{\widetilde S}+\widetilde{\Delta}+c\widetilde{C}+ \widetilde{E}=
\varphi^*(K_S+\Delta +cC),
\end{equation*}
where $\widetilde E$ is the exceptional divisor, 
$\widetilde C$ and $\widetilde \Delta$ are proper transforms of 
$C$ and $\Delta$, respectively. By
Corollary \ref{Msma}, $\Diff_{\widetilde{E}}(\widetilde{\Delta}+
c\widetilde{C})\in\Msm^{5/6}$. On the other hand, $K_{\widetilde E}+
\Diff_{\widetilde{E}}(\widetilde{\Delta}+ c\widetilde{C})\equiv 0$. 
By Lemma \ref{P1}, 
$\Diff_{\widetilde{E}}(\widetilde{\Delta}+ c\widetilde{C})\in [0, 5/6]$. 
Clearly, $\widetilde{E} \cap \widetilde{C}\neq \emptyset$. 
Applying Corollary~\ref{coeff_Diff} to our situation we obtain
$1-1/m+c/m\le 5/6$ for some $m\in\NN$. This yields $c\le 5/6$, a
contradiction.
\end{proof}

\begin{lemma}[cf. {\cite{Sh1}}]
\label{pair-discr}
Let $(S\ni o,\Lambda=\lambda_1\Lambda_1+\lambda_2\Lambda_2)$ be a
log surface germ such that $\lambda_1,\lambda_2\ge 5/6$. Assume
that $\discr(S,\Lambda)\ge -5/6$ at $o$. Then
$\lambda_1+\lambda_2\le 11/6$.
\end{lemma}
\begin{proof}
By Lemma \ref{LCT-2}, $K_S+\Lambda_1+\Lambda_2$ is lc at $o$. In
this situation there is an analytic isomorphism 
(cf. Proposition \ref{iso})
\[
(S,\Lambda, o)\simeq (\CC^2,\{xy=0\},0)/\ZZ_m(1,q),
\]
where $m\in\NN$ and $\gcd(m,q)=1$. Take $q$ so that $1\le q< m$
and consider the weighted blow up with weights $\frac1m(1,q)$. We
get the exceptional divisor $E$ with discrepancy
\[
-\frac56\le
a(E,\Lambda)=-1+\frac{1+q}m-\frac{\lambda_1}m-\frac{q\lambda_2}m.
\]
Thus
\begin{multline*}
0\le 1+q- \lambda_1-q\lambda_2-\frac m6\le 1+q 
-\frac56(1+q)-\frac m6=\frac{1+q-m}{6}.
\end{multline*}
If $m\ge 2$, this gives as $q=m-1$ and equalities
$\lambda_1=\lambda_2= 5/6$. In the case $m=1$, $q=1$ we have $0\le
2- \lambda_1-\lambda_2-1/6$, i.e. $\lambda_1+\lambda_2\le 2-1/6$.
\end{proof}

\section{Proof of the main theorem}
\label{proof}
In this section we prove Theorem \ref {main}.
First we reduce the problem to the case when $X$ 
is $\QQ$-factorial and has only log terminal singularities.
These arguments are quite standard, so the reader can skip them.

\begin{lemma}
Let $(X,\Omega)$ be a $d$-dimensional lc log variety such 
that $\Omega\in\Msm$
and let $F$ be an effective Weil $\QQ$-Cartier divisor on $X$.
Assume that the log MMP in dimension $d$ holds.
Then there is a $\QQ$-factorial $d$-dimensional klt variety  $X'$
and an effective Weil $\QQ$-Cartier divisor $F'$ on $X'$
such that  $c(X,\Omega,F)=c(X',F')$.
\end{lemma}
\begin{proof}
We prove our lemma by induction on $d$.
Put $c:=c(X,\Omega,F)$. Clearly, we may assume that 
$0<c<1$. 
Consider  minimal dlt $\QQ$-factorial 
modification $g\colon (\tilde X,\tilde \Omega)\to (X,\Omega)$ 
(see \cite[17.10]{Ut}).
By definition, this is a birational morphism $g\colon \tilde X\to X$
such that $\tilde X$ is $\QQ$-factorial and 
\[
K_{\tilde X}+\tilde \Omega+\sum E_i=g^*(K_X+\Omega)
\]     
is dlt, where $\tilde \Omega$ is the proper transform of 
$\Omega$ and the $E_i$ are prime exceptional divisors 
(if $(X,\Omega)$ is dlt, one can take $\sum E_i=0$).
Since $c>0$ and because $a(E_i,\Omega)=-1$, $F$ cannot contain 
$g(E_i)$. Therefore the proper transform of $F$
coincides with its pull-back $g^*F$. 
Replace $(X,\Omega,F)$ with $(\tilde X,\tilde \Omega,g^*F)$. 
From now on we may assume that $(X,\Omega)$ is dlt and $X$ is 
$\QQ$-factorial. There is an exceptional divisor $E$ such that 
$a(E,\Omega+cF)=-1$ and $a(E,\Omega)>-1$. 
Regard $X$ as a germ near a point
$P\in \mt{Center}(E)$.

Assume that $\down{\Omega}\neq 0$. Let $S$ be a component of 
$\down{\Omega}$ (passing through $P$).  Then $(S,\Diff_S(\Omega-S))$
is lc \cite[17.7]{Ut} and $\Diff_S(\Omega-S)\in\Msm$ 
(see Corollary~\ref{Msma}). Then it is easy to see that
$c(X,\Omega,F)=c(S,\Diff_S(\Omega-S),F\bigl|_S)$.  
Taking into account $\mathcal{T}_{d-1}\subset\mathcal{T}_d$
(see \cite[8.21]{Ko}), we get our assertion.

Now consider the case $\down{\Omega}= 0$.
Then $(X,\Omega)$ is klt. Since $X$ is a germ near $P$, 
$n(K_X+\Omega)\sim 0$ for some 
$n\in \NN$. Take $n$ to be minimal with this property. Then 
the isomorphism $\OOO_X(n(K_X+\Omega))\simeq \OOO_X$ defines 
an $\OOO_X$-algebra structure on 
$\sum _{i=0}^{n-1}\OOO_X(\down {-iK_X-i\Omega})$ this gives us 
a cyclic $\ZZ_n$-cover 
\[
\varphi\colon X':=\mt{Spec}\left(
\sum _{i=0}^{n-1}\OOO_X\left(\down {-iK_X-i\Omega}\right)\right)
\longrightarrow X. 
\]
The ramification divisor of $\varphi$ is $\Omega$.
Hence $\varphi^*(K_X+\Omega)=K_{X'}$ and $X'$ has only
log terminal singularities \cite[20.3]{Ut}. 
Put $F':=\varphi^*F$. Then 
$c(X,\Omega,F)=c(X',F')$ (see \cite[8.12]{Ko}).
Replacing $X'$ with its $\QQ$-factorialization we get 
the desired log pair.
\end{proof}

\subsection{Notation}
Let $X$ be a three-dimensional $\QQ$-factorial normal 
variety with only log
terminal singularities and let $F$ be an effective Weil
$\QQ$-Cartier divisor on $X$. Put $c:=c(F,X)$. 
Let $f\colon Y\to X$ be an inductive blowup of the pair $(X,cF)$.
Write $f^*(K_X+cF)=K_Y+cF_Y+S$, where $F_Y$ is the proper
transform of $F$ on $Y$ and $S$ is the exceptional divisor. Let
$\Theta:=\Diff_S(cF_Y)$ and $\Theta=\sum \vartheta_i\Theta_i$.

\subsection{Main assumption}
\label{epsilon}
Fix $\epsilon>0$ and assume that $1>c>5/6+\epsilon$. We prove that
there are only a finite number of possibilities for such $c$.

\begin{lemma}
\label{point}
$f(S)$ is a point.
\end{lemma}
\begin{proof}
Otherwise $f(S)$ is a curve and the pair $(X,cF)$ is lc but not
klt along $f(S)$. 
Taking a general hyperplane section we derive a
contradiction with Lemma \ref{LCT-2}.
\end{proof}

\begin{lemma}
$(Y,S+cF_Y)$ is plt.
\end{lemma}
\begin{proof}
Assume the converse. Then there is an exceptional divisor $E$ such
that $a(E,S+cF_Y)=-1$. Since $(Y,S)$ is plt,
$\mt{Center}(E)\subset E\cap F_Y$.

If  $\mt{Center}(E)$ is a curve, then $(Y,S+cF_Y)$ is lc but
not klt along $\mt{Center}(E)$. 
As in the proof of Lemma \ref{point}  
we derive a contradiction. Thus 
we may assume that $(Y,S+cF_Y)$ is plt in codimension two. 
By Adjunction \cite[Th.
17.6]{Ut} this implies that $\down{\Theta}=0$. 

Hence $\mt{Center}(E)$ is a point. Again by Adjunction
$(S,\Theta)$ is lc but not klt near $\mt{Center}(E)$. 
As above, we have a contradiction with
Lemma \ref{LCT-2}.
\end{proof}

\begin{corollary}
$(S,\Theta)$ is klt.
\end{corollary}

\subsection{}
Now we are going to construct a ``good'' birational model
$(\bar S,\bar \Theta)$ of $(S,\Theta)$. The construction is similar 
to that in \cite{Sh1}. 
Assumption \ref{epsilon} gives us that $\Theta \in \Msm^{5/6}$.
If $\discr(S,\Theta)\ge -5/6$ and $\rho(S)=1$, we put 
$(\bar S,\bar \Theta)=(S,\Theta)$. 

From now on we assume either $\discr(S,\Theta)< -5/6$ or $\rho(S)>1$,
Since $(S,\Theta)$
is klt, there is only a finite set $\EEE$ of divisors $E$ with
$a(E,\Theta)< -5/6$ \cite[2.12.2]{Ut}. Let $\mu\colon \widetilde S\to
S$ be the blow-up of all divisors $E\in \EEE$ (see \cite[Th.
17.10]{Ut}) and let $\widetilde\Theta$ be the crepant pull-back:
\[
K_{\widetilde S} +\widetilde \Theta= \mu^*(K_S+ \Theta),\quad
\mu_*\widetilde\Theta=\Theta.
\]
Then $\discr(\widetilde S,\widetilde\Theta)\ge -5/6$ and again we have
$\widetilde \Theta\in\Msm^{5/6}$. Write $\widetilde\Theta= \sum
\vartheta_i\widetilde\Theta_i$ and consider the boundary 
$\widetilde \Xi$ with $\Supp(\widetilde \Xi)=\Supp(\widetilde \Theta)$:
\[
\widetilde \Xi:=\sum \xi_i\widetilde\Theta_i,\qquad \xi_i=
\begin{cases}
1 & \text{if $\vartheta_i> 5/6$}, \\ \vartheta_i &
\text{otherwise}.
\end{cases}
\]
For sufficiently small positive $\alpha$, the $\QQ$-divisor
$\widetilde \Theta -\alpha (\widetilde \Xi -\widetilde \Theta)$
is a boundary. 
It is clear
that 
\[
K_{\widetilde S}+\widetilde \Theta -\alpha (\widetilde \Xi -
\widetilde \Theta)\equiv -\alpha (\widetilde \Xi -
\widetilde \Theta)
\]
cannot be nef. By our assumption, $\rho(\widetilde S)>1$.
Note also that $(\widetilde S,\widetilde \Xi)$ is lc
(see Lemma~\ref{LCT-2}).
Run 
$K_{\widetilde S}+\widetilde \Theta -\alpha (\widetilde \Xi -
\widetilde \Theta)$-MMP. On each step we contract an extremal ray 
$R$ such that 
\[
(K_{\widetilde S}+\widetilde \Xi)\cdot R=(\widetilde \Xi -
\widetilde \Theta)\cdot R>0.
\] 
Consider such a contraction $\varphi\colon\widetilde S\to
S^\sharp$.

\subsection{}
\label{dim=1}
Assume that $\dim S^\sharp=1$ and let $C$ be a general fiber.
Since $(\widetilde \Xi-\widetilde\Theta)\cdot C>0$, there is a component
$\widetilde\Theta_i$ of $\widetilde\Theta$ with coefficient $\vartheta_i>
5/6$ meeting $C$. 
Hence $\Diff_C(\widetilde \Theta)$ also has a component
with coefficient $>5/6$.  By Adjunction
$K_C+\Diff_C(\widetilde \Theta)$ is klt.  On the other hand,
\begin{equation*}
\label{new} 
K_{C}+\Diff_C(\widetilde \Theta)\equiv 0\quad \text{and}
\quad
\Diff_C(\widetilde \Theta)\in \Msm^{5/6}
\end{equation*}
(see Corollary~\ref{Msma}). This contradicts
Lemma \ref{P1}.

Thus, $\varphi$ is birational. 

\subsection{}
\label{dim=2}
We claim that 
$\varphi$ cannot contract a component of $\down {\widetilde
\Xi}$. Indeed, assume that $\varphi$ contracts a curve 
$C\subset \down {\widetilde
\Xi}$.  
Take $\widetilde \Theta':=\widetilde \Theta+\alpha C$
so that $\down{\widetilde \Theta'}=C$ and  
$\widetilde \Theta'\le \widetilde \Xi$.
Since $C^2<0$, we have 
$(K_{\widetilde S}+\widetilde \Theta')\cdot C<0$.
Therefore 
\[
\left(K_{\widetilde S}+\widetilde \Theta'+
\beta \bigl(\widetilde \Xi-\widetilde \Theta'\bigr)\right)
\cdot C=0
\]
for some $0<\beta<1$. 
Put $\widetilde \Theta'':=\widetilde \Theta'+
\beta (\widetilde \Xi-\widetilde \Theta')$.
Then $\widetilde \Theta''\le \widetilde \Xi$, so 
$(\widetilde S,\widetilde \Theta'')$ is lc.
Moreover $\widetilde \Theta''\in\Msm^{5/6}$.
Since $(\widetilde \Xi-\widetilde\Theta'')\cdot C>0$, 
there is a component of $\widetilde \Xi-\widetilde\Theta''$
meeting $C$. By Lemma \ref{LCT-2}, 
$(\widetilde S,\widetilde\Theta'')$ is plt
near $C\cap \Supp(\widetilde \Xi-\widetilde\Theta'')$.
As in \ref{dim=1} we derive a contradiction by Lemma \ref{P1}.

Put $\Xi^\sharp:=\varphi_{*}\widetilde\Xi$ and
$\Theta^\sharp:=\varphi_{*}\widetilde\Theta$. By \cite[2.28]{Ut},
\[
\discr(S^\sharp,\Theta^\sharp)=\discr(\widetilde S,\widetilde \Theta) \ge
-5/6.
\]
Thus all the assumptions hold for $(S^\sharp,\Theta^\sharp)$. Again
\[
K_{S^\sharp}+\Theta^\sharp -\alpha (\Xi^\sharp -
\Theta^\sharp)\equiv -\alpha (\Xi^\sharp -
\Theta^\sharp)
\]
cannot be nef.

Continuing the process we get a new pair $(\bar S, \bar
\Theta)$ such that
\par\medskip\noindent
$\rho(\bar S)=1$, $\bar \Theta\in \Msm^{5/6}$, 
$(\bar S,\bar \Theta)$ is klt,
$K_{\bar S}+\bar \Theta\equiv 0$, and
$\discr(\bar S,\bar \Theta)\ge -5/6$.
\par\medskip\noindent
Note that all our birational modifications are 
$(K+\Theta)$-crepant.
Hence
\[
\totaldiscr(S,\Theta)=\totaldiscr(\bar S,\bar \Theta)
=\totaldiscr(\widetilde S, \widetilde \Theta)
\]
(see \cite[3.10]{Ko}).
Consider the decomposition $\Theta=\Theta^a+\Theta^b$, where
\[
\Theta^a=\sum_{\Theta_i\subset F_Y} \vartheta_i\Theta_i, \qquad 
\Theta^b=\sum_{\Theta_i\not \subset F_Y} \vartheta_i\Theta_i.
\]
Similarly, $\bar \Theta=\bar\Theta^a+\bar\Theta^b+\bar\Theta^c$, 
where $\bar\Theta^a$ and $\bar\Theta^b$ are proper transforms 
of $\Theta^a$ and $\Theta^b$, respectively, and components of
$\bar\Theta^c=\bar \Theta-\bar\Theta^a-\bar\Theta^b$ are 
proper transforms of exceptional divisors of $\mu$.

It is clear that $\Theta^b, \bar\Theta^b\in \Msm$
and $\bar\Theta^c\in (5/6,1)$. 
Since the coefficients of $\Theta^a$ (as well as $\bar\Theta^a$)
are of the form 
\[
\vartheta_i=1-1/m_i+c/m_i\ge c>5/6+\epsilon,
\]
we have $\Theta^a\in (5/6+\epsilon,1)$.
By our assumptions $\Theta^a\neq 0$.

We need the following result of Alexeev \cite{A}:
\begin{theorem}
\label{Alexeev}
Fix $\epsilon>0$.
Consider the class of all projective log surfaces $(S,\Theta)$
such that $-(K_S+\Theta)$ is nef and 
$\totaldiscr(S,\Theta)>-1+\epsilon$ excluding only the case
\begin{itemize}
\item
$\Theta=0$, $K_S\equiv 0$ and the singularities of $S$ are at 
worst Du Val.
\end{itemize}
Then the class $\{S\}$ is bounded, i.e. $S$ belongs to a finite 
number of algebraic families. 
\end{theorem}

\subsubsection{}
\label{prep}
Let $\bar \Theta_1$  be a component of  $\bar \Theta^a$.
Then $\vartheta_1>5/6+\epsilon$. 
Since $\rho(\bar S)=1$, every two components of $\bar \Theta$
intersects each other. Applying Lemma \ref{pair-discr} we 
obtain
\[
\vartheta_j\le 11/6-\vartheta_1<11/6-5/6-\epsilon=1-\epsilon
\]
for all $j\neq 1$.
Since $\bar \Theta^b\in \Msm$, there is  
only a finite number of possibilities for the coefficients 
of $\bar \Theta^b$ (and $\Theta^b$). 

\subsubsection{}
\label{not_irred}
If $\bar \Theta^a$ has  at least two components,
say $\bar \Theta_1$ and  $\bar \Theta_2$, then 
by Lemma~\ref{pair-discr} 
the inequality $\vartheta_k<1-\epsilon$ holds for
all $\vartheta_k$. Thus 
\[
\totaldiscr(S,\Theta)=\totaldiscr(\bar S,\bar 
\Theta)>-1+\epsilon.
\] 
Apply \ref{Alexeev} to $(S,\Theta)$. 

For all coefficients of $\Theta$ we have 
$\vartheta_i\ge 1/2$. Fix a very ample divisor $H$ on $S$.
Then $H\cdot \sum \Theta_i\le 2 H\cdot K_S\le \mt{Const}$. 
This shows that the pair $(S,\Supp(\Theta))$ is also bounded.

As above, 
$(S,\Supp(\Theta))$ is bounded. From the equality
$0=K_S^2+K_S\cdot \Theta^a+K_S\cdot \Theta^b$ we obtain 
\[
\sum_{\Theta_i\not \subset F_Y} (1-1/m_i+c/m_i) (K_S\cdot\Theta_i)=
-K_S^2-K_S\cdot \Theta^b,
\]
where $1-1/m_i+c/m_i<1-\epsilon$.
This gives us a finite number of possibilities for
$c$.

\subsubsection{}
Assume that $\bar \Theta^a=\vartheta_1\bar\Theta_1$, 
where $\vartheta_1=1-1/m_1+c/m_1$. If $\vartheta_1<1-\epsilon$,
then we can argue as above. Let $\vartheta_1\ge 1-\epsilon$.
Then $\Theta_1$ is the only
divisor with discrepancy $a(\Theta_1, \Theta)\le -1+\epsilon$.
Put $\Lambda:=\Theta -\vartheta_1\bar\Theta_1$. 
Then $a(\Theta_1, \Lambda)=0$, so
$\totaldiscr(S,\Lambda)>-1+\epsilon$.
Note that $\Theta_1$ is ample (because $\Theta_1=(F_Y\bigl|_S)_{\red}$
and $F_Y$ is $f$-ample, see \ref{constr-plt-blowup}, (iii)).
Hence $-(K_S+\Lambda)$ is also ample. By \ref{Alexeev}
$(S,\Supp(\Lambda))$ is bounded and so is $(S,\Supp(\Theta))$. 
As in \ref{not_irred}, there is only a finite number of possibilities for
$c$.

The following example illustrates our proof:

\begin{example}
\label{ex-2}
Notation as in Example \ref{ex-1}. Assume that $\gcd(6,r)=1$. Let
$f\colon Y\to X$ be the weighted blowup with weights $(3r,2r,6)$.
Then $f$ is an inductive blowup of $(X,cF)$ and the exceptional
divisor $S$ is isomorphic to $\PP(3r,2r,6)\simeq\PP^2$. It is easy
to compute that
$\Theta=\Diff_S(cF_Y)=\frac12L_1+\frac23L_2+\frac{r-1}rL_3+cL_0$, where
$c=5/6+1/r$ and 
$L_1,L_2,L_3,L_0$ are lines on $S\simeq\PP^2$ given by equations
$x=0$, $y=0$, $z=0$ and $x+y+z=0$, respectively. Thus 
$\discr(S,\Theta)\ge -5/6$ and $\bar S=\widetilde S=S\simeq \PP^2$.
\end{example}

\subsection*{Concluding remark} 
(i)
Using the same arguments one can see that see 
that the set $\mathcal{T}_3$ in  Theorem \ref{main}
can be replaced with  $\mathcal{T}_3(\Msm)$, the set of all 
values $c(X,\Omega,F)$ with $\Omega\in\Msm$.

(ii) 
We expect that 
our proof of Theorem \ref{main} can be generalized in higher dimensions
modulo the following facts: the log MMP, boundedness result \ref{Alexeev}
and lemmas \ref{P1} and \ref{LCT-2}.
Also we hope that our method allow us to get the complete 
description of $\mathcal{T}_3\cap [5/6,1]$.

\end{document}